\documentclass[11pt]{article}
\usepackage{amsthm}
\usepackage{thmtools}

\makeatletter
\def\@seccntformat#1{%
  \expandafter\ifx\csname c@#1\endcsname\c@section\else
  \csname the#1\endcsname\quad
  \fi}
\makeatother

\newtheorem{theorem}{Theorem}[section]

\newtheorem{definition}[theorem]{Definition}

\theoremstyle{definition}

\newcommand{\DShortExactSeq}[3]{\xymatrix@1{0\ar[r]&{{#1}}\ar[r]&{{#2}}\ar[r]&{{#3}}\ar[r]&0}}

\usepackage{xcolor}
\definecolor{darkgreen}{RGB}{34,139,34}

\usepackage{mathrsfs}
\usepackage{bbm}

\usepackage{adjustbox}
\usepackage{times}
\usepackage{tikz}
\usepackage{verbatim}
\usepackage{mathtools}
\usepackage{fullpage}
\usepackage{graphicx} 
\usepackage{pgfplots} 
\usepackage{fancyhdr} 
\usepackage{lastpage} 
\usepackage{colortbl}
\usepackage{ulem}
\usepackage{hyperref}
\usepackage{thmtools} 
\usepackage{wrapfig}
\usepackage{pifont} 
\usepackage{scrextend} 
\usepackage{pgfplotstable}
\usepackage{amssymb, amsmath}
\usepackage{polynom}
\usepackage{mathrsfs}
\usepackage{makeidx}
\usepackage{tabularx}
\usepackage{algorithm2e}
\usetikzlibrary{patterns}
\usetikzlibrary{shapes,snakes}
\usetikzlibrary{3d}

\makeindex
\title{Threaded Gr\"{o}bner Bases: a {\it Macaulay2} package {\it ThreadedGB}}
\author{Sonja Petrovi\'c and Shahrzad Jamshidi Zelenberg}

\begin{document}
\maketitle

\begin{abstract} 
      The complexity of Gr\"{o}bner  computations has inspired many improvements to Buchberger's 
      algorithm over the years.    Looking for further insights into the algorithm's performance, we offer a threaded implementation of classical Buchberger's algorithm in {\it Macaulay2}. 
The output of the main function of the package includes information about  {\it lineages}  of non-zero remainders that are added to the basis during the computation. 
This information can be used  for further algorithm improvements and optimization. 
  \end{abstract}

\section{Introduction}

The importance in computational algebra of Gr\"obner bases and therefore of Buchberger's algorithm, as well as its many variants, is indisputable. Yet it is still a challenge to apply brute force algorithms to larger problems primarily due to considerations in computer science. That is, the current computing paradigm favors clusters of CPUs, or nodes, rather than one massive CPU. As a result, there is a need to distribute this algorithm that is automated for the user (in that it does not require a user to know how it should be distributed). 

Past work in this area has focused on synchronized methods as detailed in \cite{mityunin2007}. One method spreads a key step in Buchberger's algorithm---the reduction of S-pairs by division---across nodes; another sends tasks to individual nodes while a central, coordinating node waits for all threads to complete. Each of these still requires some central node and synchronization, which leads to bottlenecks in the computation. A truly distributed algorithm would be decentralized and asynchronous. 
Within \cite{SJZ_thesis},
Zelenberg discusses an asynchronous, decentralized distributed version of Buchberger's algorithm done generically with the potential of very good speedups. Zelenberg implemented a threaded version  in {{\it Python}} \cite{python} 
to explore this further,  and as a result some important discoveries were made. 
It should be noted  that  multithreaded algorithms are not necessarily distributed across distinct nodes; rather, threads are sharing computation and passing information back and forth.

Of the discoveries made, the most important is this: in order for a distributed process to be both generically usable and automated for the user, an effective algorithm will need to account for features of the polynomials (relative to the starting basis) when deciding what tasks to assign to what node. This is because transferring information between nodes is a very slow process and therefore needs to be minimized. 

Given the need to analyze aspects of these polynomials, {\it Macaulay2} \cite{M2} offers some clear advantages over {\it Python}. Moreover, the {\it Macaulay2}  engine is written in C/C++, a language well suited for writing distributed algorithms. 
Threads within {\it Macaulay2} work differently than within {\it Python} and, as such, some design changes were necessary from the implementation discussed in \cite{SJZ_thesis}. Queues are replaced altogether with  
 hash tables---an improvement since threads access the most up-to-date version of the generating set at the time of reduction. Even so, this cannot eliminate redundancies as threads may compute the same result (virtually) simultaneously. 

One of the goals of our package {\it ThreadedGB} is to allow a user to analyze what we refer to as {\it lineages} of polynomials in a Gr\"obner basis. 
\begin{definition} 
Let $G$ be a  Gr\"obner basis of $I=(f_0,\dots,f_k)$.  
 A \emph{lineage} of a polynomial in $G$  is a natural number, or an ordered pair of lineages, tracing      its history in the given Gr\"obner basis computation. It is defined recursively as follows: 
 \begin{itemize}
 \item  For the starting generating set,  $Lineage(f_i) = i$, 
 \item For any subsequently created S-polynomial $S(f,g)$, the lineage of its remainder $r$ on division is the pair 
 $Lineage(r)=(Lineage(f),Lineage(g))$. 
 \end{itemize}
 
\end{definition}

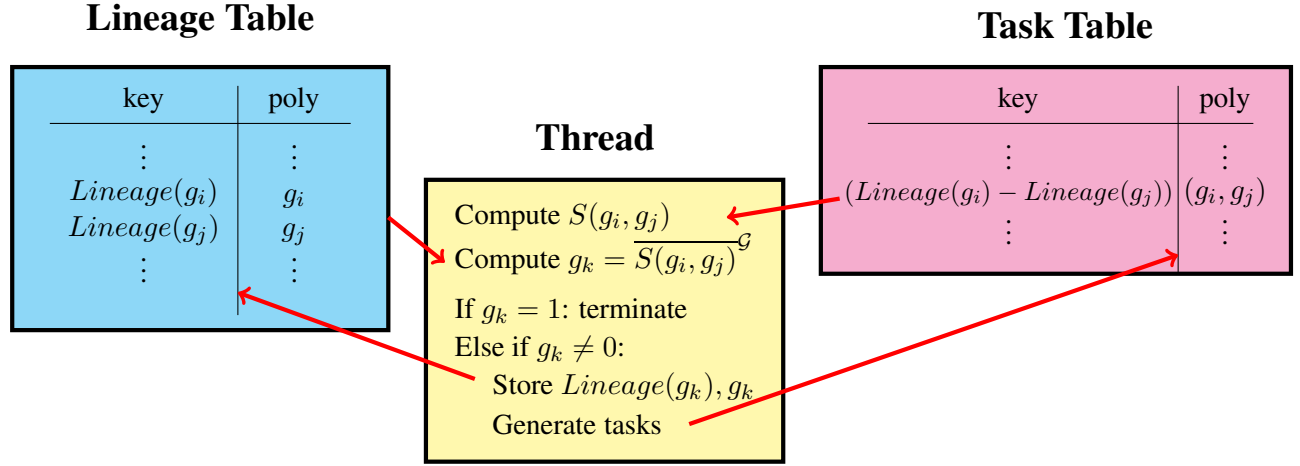
\begin{figure}
\begin{center}\begin{tikzpicture}
\node[above] at (1.5, 3.75){\Large\bf Lineage Table}; 
\draw[ultra thick, fill = cyan!40] (-1, 3.5) rectangle (4, 0);
\draw (-.5, 2.75) -- (3.5, 2.75);
\draw (2, 3.25) -- (2, 0.2);
\node[above] at (.75, 2.75){key}; 
\node[above] at (2.75, 2.75){poly}; 
\node[above] at (.75, 2){$\vdots$}; 
\node[above] at (2.75, 2){$\vdots$}; 
\node[above] at (.75, 1.5){$Lineage(g_i)$}; 
\node[above] at (2.75, 1.5){$g_i$}; 
\node[above] at (.75, 1){$Lineage(g_j)$}; 
\node[above] at (2.75, 1){$g_j$}; 
\node[above] at (.75, 0.5){$\vdots$}; 
\node[above] at (2.75, 0.5){$\vdots$};


\node[above] at (6.75, 2.25){\Large\bf Thread}; 
\draw[ultra thick, fill = yellow!40] (4.5, 2) rectangle (9.25, -1.75);
\node[above, right] at (4.75, 1.5){Compute $S(g_i, g_j)$}; 
\node[above, right] at (4.75, 1){Compute $g_k = \overline{S(g_i, g_j)}^\mathcal{G}$}; 
\node[above, right] at (4.75, 0.25){If $g_k = 1$: terminate}; 
\node[above, right] at (4.75, -0.25){Else if $g_k\neq 0$:};
\node[above, right] at (5.25, -0.75) {Store $Lineage(g_k), g_k$}; 
\node[above, right] at (5.25, -1.25){ Generate tasks};


\node[above] at (13, 3.75){\Large\bf Task Table}; 
\draw[ultra thick, fill = magenta!40] (9.75, 3.5) rectangle (16, 0.75);
\draw (10.375, 2.75) -- (15.5, 2.75);
\draw (14.5, 3.25) -- (14.5, 0.75);
\node[above] at (12.375, 2.75){key}; 
\node[above] at (15.125, 2.75){poly}; 
\node[above] at (12.25, 2){$\vdots$}; 
\node[above] at (15.125, 2){$\vdots$}; 
\node[above] at (12.25, 1.5){\small $(Lineage(g_i) - Lineage(g_j))$}; 
\node[above] at (15.125, 1.5){$(g_i, g_j)$}; 
\node[above] at (12.25, 1){$\vdots$}; 
\node[above] at (15.125, 1){$\vdots$};

\draw[ultra thick, red, ->] (10, 1.75) -- (8.5, 1.5);
\draw[ultra thick, red, ->] (4, 1.5) -- (4.75, .9);
\draw[ultra thick, red, ->] (5.15, -0.65) -- (2, 0.5);
\draw[ultra thick, red, ->] (8, -1.25) -- (14.5, 1);

\end{tikzpicture}\end{center}
  \caption{
  The lineage table is initialized with the starting basis, and the task
     table is initialized with a task for each pair of starting generators. 
  Then, all threads perform the  same task: they pull a pair of polynomials from the tasks hash table along with an associated lineage key. They compute the S-polynomial and reduce it with respect to the current basis. If the remainder, $r$, is nonzero, it is stored in  the lineage table along with its lineage key, and, for each $g$ in the current basis, a new task indexed by the pair $(g,r)$ is added to the task table. 
    If the remainder 1 is found, the process of creating tasks stops. The process is repeated using $n$ parallel threads, where $n$ is specified by the user, until the task table is empty.  
}\label{fig:tgb}
\end{figure}

To illustrate, suppose $I=(x^2 - y, x^3 - z)\subset \mathbb Q[x,y,z]$ with graded reverse lexicographic order. Then $Lineage(x^2-y)=0$ and $Lineage(x^3-z)=1$.   
 Two additional elements are  added to create a (non-minimal) Gr\"obner basis: $xy + z$ and $y^2 -xz$, with lineages $((0,1),0)$ and $(0,1)$, respectively.  According to $Lineage(y^2-xz)$, this element is constructed from  $S(xy + z, x^2 - y)$. 
 Lineages are expressions of the starting basis and thus dependent on the choice and order of its elements. 
More importantly, a lineage is not necessarily unique, as the same polynomial can be constructed multiple ways. The lineage tables produced by {\it ThreadedGB} (see Figure~\ref{fig:tgb}) do not provide all possible lineages---only a particular choice based on the order in which the basis elements are provided by the user.

\section{Effect of ordering of polynomials on lineages: a simple example}

Consider $\mathbb Q[x_1,x_3,x_0,x_4,x_2]$ with lexicographic order and the ideal of the rational normal curve in $\mathbb P^3$.
The six generators have lineages $0\dots 5$, and Buchberger's algorithm adds 3 new elements  to the Gr\"obner basis before final reduction. This can be seen by turning on the {\tt gbTrace} option in {\it Macaulay2}, which tells us three new polynomials are added to the basis. The function {\tt tgb} lets us know exactly which ones and their lineages. Specifically, a run of {\tt tgb} reveals these are $ x_0x_4-x_2^2, -x_3x_0x_4+x_3x_2^2,
      -x_0x_4x_2+x_2^3
$
with lineages  $(2,3), (1,4), (1,2)$, respectively. 
\small
\begin{verbatim}
i4 : QQ[x_1,x_3,x_0,x_4,x_2,MonomialOrder=>Lex];
i5 : rnc = minors(2, matrix{{x_0..x_3},{x_1..x_4}});
o5 : Ideal of QQ[x , x , x , x , x ]
                  1   3   0   4   2
i6 : allowableThreads = 4; 
i7:  g = tgb(rnc)
                                        3
o7 = LineageTable{(1, 2) => - x x x  + x   }
                               0 4 2    2
                                          2
                  (1, 4) => - x x x  + x x
                               3 0 4    3 2
                                    2
                  (2, 3) => x x  - x
                             0 4    2
                          2
                  0 => - x  + x x
                          1    0 2
                  1 => - x x  + x x
                          1 2    3 0
                               2
                  2 => x x  - x
                        1 3    2
                  3 => - x x  + x x
                          1 3    0 4
                  4 => x x  - x x
                        1 4    3 2
                          2
                  5 => - x  + x x
                          3    4 2
o7 : LineageTable
\end{verbatim}
\normalsize
Running the command {\tt reduce g} will produce a reduced Gr\"obner basis; in particular, the lineage table entries with keys $(1-2)$, $(1-4)$ and $2$ will be replaced by {\tt null}.  
This allows the user to see which non-zero polynomials produced during the computation turn out not to be needed.
Of course, to continue computing with the given basis, one wishes to have it in standard {\it Macaulay2} format, which is a matrix. 
\small
\begin{verbatim}
i7 : matrix reduce g
o7 = | x_1^2-x_0x_2 x_1x_2-x_3x_0 x_1x_3-x_2^2 
        x_1x_4-x_3x_2 x_3^2-x_4x_2 x_0x_4-x_2^2 |
                                     1                              6
o7 : Matrix (QQ[x , x , x , x , x ])  <--- (QQ[x , x , x , x , x ])
                  1   3   0   4   2              1   3   0   4   2
\end{verbatim}
\normalsize

One can use the package to study, for example,  how reordering the input basis affects the algorithm. 
In Gr\"obner computations, {\it Macaulay2} creates and processes S-polynomials in lexicographic order of pairs (first and second, then first and third, and so on). 
Let $S=\mathbb Q[a,b,c,d]$ and $I=( abc-1,abc, -c^3+a^2+bd)$; clearly $I=S$. But the order of generators listed affects the complexity of the particular run; namely, listing the cubic first makes the algorithm perform more steps.
The method {\tt tgb} can be verbose and can tell us what is going on behind the scenes for each lineage.  

\small
\begin{verbatim}
i8 : QQ[a, b, c, d]; 
i9 : I =  ideal (a*b*c,a*b*c - 1,  - c  + a  + b*d); 
i10 :  tgb(I,Verbose=>true)
Scheduling a task for lineage (0,1)
Scheduling a task for lineage (0,2)
Scheduling a task for lineage (1,2)
Adding the following remainder to GB: 1 from lineage (0,1)
Found a unit in the Groebner basis; reducing now.
o10 = LineageTable{(0, 1) => 1}
                  0 => null
                  1 => null
                  2 => null
o10 : LineageTable
\end{verbatim}
\normalsize
Compare this to the following run of threaded Buchberger's algorithm under a different input order.  
\small
\begin{verbatim}
i11 : I = ideal (-c+a+b*d, a*b*c-1, a*b*c);
o11 : Ideal of QQ[a, b, c, d]
i12 : tgb(I) 
o12 = LineageTable{(0, 1) => null}
                  (0, 2) => null
                  (1, 2) => 1
                  0 => null
                  1 => null
                  2 => null
o12 : LineageTable
\end{verbatim}
\normalsize
Three new elements are added to the basis, namely (0,1), (0,2), (1,2), if the cubic generator is listed first, but if it is listed last, then only the polynomial with lineage (0,1) is added---because it already equals $1$---and the algorithm stops.

\section{Nuts and Bolts}

      Given a list $L$ or an ideal $I$ and an integer $n$, the main method {\tt tgb} 
       uses  Tasks in {\it Macaulay2} to compute a Gr\"obner basis of $I$ or $(L)$ 
      using $n$ threads. 
It returns an object of type LineageTable, which is an instance of HashTable, whose values are a Gr\"obner basis of $I$ or $(L)$. 
      The keys are polynomial lineages. 
      
      The starting basis $L$ (meaning, the input list {\tt L} or {\tt L=gens I}) 
      populates 
      the entries of a lineage table  $G$ with keys from $0$ to one less than the number of elements of $L$.   
       The method creates all possible
      S-polynomials of $L$ and schedules their reduction with respect to $G$ as tasks. 
      Throughout the computation, every nonzero remainder added to the basis is added to $G$ 
      with its lineage, as defined above, being the key. Each such remainder also triggers the creation of S-polynomials
      using it and every element in $G$ and scheduling the reduction thereof as additional tasks. 
      The process is done when there are no remaining tasks.
      
      There is a way to track the tasks being created by turning on the option {\tt Verbose}, or provide the reduced or a minimal Gr\"obner basis using the functions {\tt reduce} or {\tt minimalize}, respectively. The users who 
       expect just a Gr\"obner basis in usual {\it Macaulay2} format, without the lineages, can call {\tt matrix LineageTable}.

\section{Improvements and speed-ups}

As with any {\it Macaulay2} package, improvements are easy to make via GitHub. Our package's GitHub repository will be made public shortly, so other users can implement any extensions or add improvements to this threaded implementation of Buchberger's algorithm. These may include known speed-ups as optional ways to run the algorithm; for example, if one wishes to study lineages produced by the F4 algorithm \cite{faugereF4}, then one can build that option into this  threaded computation. 

The current goal is to explore algorithm performance and complexity and how input basis features affect these; the lineages are designed specifically to aid in this goal. 
Of course, speed-ups should come `naturally' from a threaded implementation; to achieve effective speed-ups in practice,   we plan to implement {\tt tgb} in the engine, using C/C++.

\small 

\nocite{misc}
\bibliographystyle{unsrt}
\bibliography{Biblio}
\printindex
\end{document}